\documentclass[preprint,12pt]{elsarticle}
\usepackage{mathrsfs}
\usepackage{amsfonts}
\usepackage{amsmath,amssymb}
\usepackage[all]{xy}
\usepackage{latexsym}
\usepackage{amsthm,a4wide,color}
\usepackage{amsmath,amscd,verbatim,enumerate}
\usepackage{hyperref}
\usepackage{graphicx}
\usepackage{extarrows}
\usepackage{lineno}
\biboptions{sort&compress}

\input diagxy

\theoremstyle{plain}
  \newtheorem{thm}{Theorem}[section]
   \newtheorem{theorem}[thm]{Theorem}
  \newtheorem{lemma}[thm]{Lemma}
  \newtheorem{proposition}[thm]{Proposition}
  \newtheorem{corollary}[thm]{Corollary}
\theoremstyle{definition}
  \newtheorem{definition}[thm]{Definition}
  
  \newtheorem{example}[thm]{Example}
  \newtheorem{remark}[thm]{Remark}

\begin{document}

\newcommand{\oto}{{\to\hspace*{-3.1ex}{\circ}\hspace*{1.9ex}}}
\newcommand{\lam}{\lambda}

\newcommand{\TC}{\mathbb{C}}
\newcommand{\TD}{\mathbb{D}}
\newcommand{\TB}{\mathbb{B}}
\newcommand{\TS}{\mathbb{S}}

\numberwithin{equation}{section}
\renewcommand{\theequation}{\thesection.\arabic{equation}}

\begin{frontmatter}
\title{Stratified $L$-convex groups\tnoteref{F}}

\tnotetext[F]{This work is supported by National Natural Science Foundation of China (No. 12171220)}

\author{Lingqiang Li\corref{cor}}
\ead{lilingqiang0614@126.com}\cortext[cor]{Corresponding author,
Tel: +86 15206506635}
\author{Qiu Jin}
\ead{jinqiu79@126.com}
\address{Department of Mathematics, Liaocheng University, Liaocheng 252059, P.R.China}
\begin{abstract} This paper investigates a novel structure of stratified $L$-convex groups, defined as groups possessing stratified $L$-convex structures, in which the group operations are $L$-convexity-preserving mappings. It is verified that stratified $L$-convex groups serve as objects, while $L$-convexity-preserving group homomorphisms serve as morphisms, together forming a concrete category, denoted as {\bf SLCG}. As a specific instance of {\bf SLCG} (i.e., when $L$=2),  the category of convex groups, denoted as {\bf CG}, is also defined. We show that  {\bf CG} can be embedded within {\bf SLCG} as a reflective subcategory. In addition, we demonstrate that {\bf SLCG} possesses well-defined characterizations, localization properties, and initial and final structures, establishing it as a topological category over groups.
\end{abstract}

\begin{keyword} $L$-convex space \sep  Stratified $L$-convex group  \sep $L$-convex remote neighborhood

{\it 2020MSC:} 52A01, 54A40, 54A20
\end{keyword}

\end{frontmatter}

\section{Introduction}

The convex structure (also known as an algebraic closure system) is an important mathematical concept \cite{MV93}, abstracted from the three basic properties of convex sets. As natural extensions, many types of (lattice-valued) fuzzy convex structures have been proposed by various scholars, including Maruyama \cite{YM09}, Pang \cite{PB23}, Shi and Xiu \cite{SFG14,SFG17}, Li \cite{LLQ17}, Sayed \cite{OS19}, and Wu \cite{WXY23}. The axiomatic definition of convex structures closely resembles that of topology. Consequently, it is quite natural to apply (fuzzy) topological methods to study (fuzzy) convex structures. Currently, lattice-valued convex structures  have received considerable attention using topological methods, including base and subbase \cite{PB19,L.L16}, closure operators \cite{SC30,HZ242}, remote neighborhood systems \cite{XZY202,HZ24}, convergence structures \cite{HXC24, WXY24,ZL20}, characterizations \cite{SS22,WXY23,HL20} and categories \cite{MY23,PB232,CX22}, among others \cite{LW24,GJ24}. The success of these results demonstrates that using lattice-valued topological methods to study lattice-valued convex structures is both convenient and practical.

A topological group is a group with a topology where the group operations are continuous \cite{AA08}.  Lattice-valued topological groups are natural extensions of topological groups, where the underlying topological space is replaced by a lattice-valued topological space. Due to the diversity of lattice-valued topologies, there exist various types of lattice-valued topological groups, such as  Bayoumi's $L$-topological group \cite{FB05}, Zhang's $L$-fuzzifying topological group \cite{Yan12},  Zhao's $(L,M)$-fuzzy topological group \cite{DZ14}, Mufarrij's  stratified $L$-topological group \cite{JA08}, Jin's Stratified lattice-valued neighborhood tower group \cite{L.L171}, Ahsanullah's quantale-valued approach group \cite{TA19}, and Li's $\top$-neighborhood group \cite{LL21}.

Motivated by the similarity between convex structures and topology, this paper primarily develops a novel structure of (lattice-valued) convex groups. Specifically, for a completely distributive De Morgan algebra $L$, we define stratified $L$-convex groups as groups with stratified $L$-convex structures, ensuring the group operations are $L$-convexity-preserving. Notably, when $L=2$, a stratified $L$-convex group reduces to a convex group. It should be noted that although lattice-valued convex groups are inspired by the definition of topological groups, there are significant differences in the form of the related results and the research methods when compared to topological groups. For more details, please refer to the subsequent Remarks \ref{remrn}, \ref{remfinal} and \ref{remcate}.


The contents are organized as follows. Section 2 recalls some basic concepts as preliminaries. Section 3 presents a further study of stratified
 $L$-convex spaces. Section 4 introduces the category of stratified $L$-convex groups (abbreviated as {\bf SLCG}), along with the category of convex groups (abbreviated as {\bf CG}). The properties and descriptions of the proposed category are examined. Section 5 examines the initial and final structures of {\bf SLCG} and demonstrates that {\bf SLCG} is a topological category over groups. Section 6 explores the relationship between {\bf CG} and {\bf SLCG}, proving that {\bf CG} can be reflectively embedded into {\bf SLCG}. Section 7 concludes the paper.

\section{Preliminaries}

Suppose $X$ is a non-empty set, and $2^X$ is its power set. The family $\{A_i\}_{i\in I}\subseteq 2^X$ is referred  directed if $\forall A_{i_1},  A_{i_2}$, there exists $ A_{i_3}$ s.t.  $  A_{i_1}\subseteq A_{i_3}$ and $ A_{i_2}\subseteq A_{i_3} $. We use  $\{A_i\}_{i\in I}{\subseteq}^\uparrow 2^X$ to denote a  directed subfamily of $2^X$, and ${\bigcup\limits_{i\in
I}}^\uparrow A_i$ to denote its join.

\begin{definition} \label{defcc} (\cite{MV93}) A subfamily $\mathcal{C}\subseteq 2^X$ is termed a convex structure on  $X$ whenever:

(1) $\emptyset, X\in \mathcal{C}$,

(2)  $\{C_k\}_{k\in K}\subseteq \mathcal{C}$ implies $\bigcap\limits_{k\in K} C_k\in \mathcal{C}$,

(3)   $\{ C_k\}_{k\in K}\subseteq^\uparrow \mathcal{C}$ implies ${\bigcup\limits_{k\in
K}}^\uparrow C_k\in \mathbb{C}$.

The pair $(X, \mathcal{C})$ is termed a convex space.
\end{definition}

A mapping  $(X,\mathcal{C})\stackrel{f}\longrightarrow (Y,\mathcal{G})$ between convex space is termed convexity-preserving (abbreviated as CP) whenever each $ G\in \mathcal{G}$ ensures $f^\leftarrow(G)\in \mathcal{C}$.

The category of convex spaces along with
CPs is denoted as {\bf CVS}.

A mapping  $(X,\mathcal{C})\stackrel{f}\longrightarrow (Y,\mathcal{G})$ between convex space is termed convex-to-convex whenever  each $ C\in \mathcal{C}$ ensures $f^\rightarrow(C)\in \mathcal{G}$.

The article assumes that $L=(L,\leq, \top, \bot, \wedge, \vee, ')$ is a completely distributive De Morgan algebra, where $\top$ and $\bot$ represent the greatest and least elements, and $'$ denotes a reversed involution negation.

An element $a\in L$ is termed coprime provided $a\leq b\vee c$ ensures $a\leq b$ or $a\leq c$. One uses $J(L)$ to represent all nonzero coprime elements in $L$ \cite{G03}.

%

An $L$-fuzzy set on $X$ refers to a mapping $A:X\longrightarrow L$. The collection of all $L$-fuzzy sets on $X$ is denoted by $L^X$. The operations
$\vee,\wedge,'$ on $L$ can be extended to $L^X$ pointwise.


For $A\in 2^X$, we use  $\top_A$ to denote its characteristic mapping, i.e., $\top_A(z)=\top$ if $z\in A$ and $\top_A(z)=\bot$ otherwise. When
$L=\{\bot, \top\}$, then $L^X$ and $2^X$  are in one-to-one correspondence via characteristic mappings.

For $a\in L$, let $\hat{a}\in L^X$ represent the constant mapping. For $x\in X$ and $a\in J(L)$, define $x_a\in L^X$ by $x_a(z)=a$ if $z=x$, and $x_a(z)=\bot$ otherwise.  The collection of all such $x_a$  is denoted by $J(L^X)$.


For $A, B\in L^X$, define the product $ A\times  B\in L^{X\times X}$ by $(x_1,x_2)\mapsto A(x_1)\wedge  B(x_2)$.

Given a mapping $f:X\longrightarrow Y$, one define
$f^{\rightarrow}: L^X\longrightarrow L^Y$ and $f^{\leftarrow}:
L^Y\longrightarrow L^X$ as
$$f^{\rightarrow}(B)(y)=\bigvee_{f(x)=y} B(x),\ \ f^{\leftarrow}(A)=A\circ f.$$

\begin{proposition} \cite{U.H99} Given a mapping $f:X\longrightarrow Y$. Then

{\rm (1)} $\forall \{C_k\}_{k\in K}\subseteq L^X$, $f^\rightarrow (\bigvee\limits_{k}C_k)=\bigvee\limits_{k}f^\rightarrow(C_k)$.

{\rm (2)} $\forall \{D_k\}_{k\in K}\subseteq L^Y$, $f^\leftarrow (\bigvee\limits_{k}D_k)=\bigvee\limits_{k}f^\leftarrow(D_k)$ and $f^\leftarrow (\bigwedge\limits_{k}D_k)=\bigwedge\limits_{k}f^\leftarrow(D_k)$.

{\rm (3)} $\forall D\in L^Y$, $f^\leftarrow (D')=[f^\leftarrow(D)]'$.

{\rm (4)} $\forall C\in L^X, D\in L^Y$, $f^\rightarrow(C)\leq D\Longleftrightarrow C\leq f^\leftarrow (D)$.

{\rm (5)} $\forall D\in L^Y$, $f^\rightarrow f^\leftarrow (D)\leq D$, and $f^\rightarrow f^\leftarrow (D)=D$ when $f$ is a surjective mapping.

\end{proposition}



A group
homomorphism is a mapping  $f:(X,\cdot)\longrightarrow (Y,\cdot)$ fulfills  $\forall u,v\in X$,
$$f(uv)=f(u)f(v),\ f(u^{-1})=f(u)^{-1}.$$

The category of group along with group homomorphism is denoted  as {\bf Grp}.

Given $x\in X$, let $$\mathcal
{L}_x:X\longrightarrow X, u\mapsto xu; \ ({\rm resp.,} \ \mathcal {R}_x:
X\longrightarrow X, u\mapsto ux),$$  denote the left (resp., right)
translation.

 For any $A, B\in L^X$, we define
$ A^{-1},  A\odot B\in L^X$ as
$$\forall u\in X,  A^{-1}(u)= A(u^{-1}), \ \ ( A\odot B)(u)=\bigvee_{xy=u}\big( A(x)\wedge  B(y)\big).$$

\begin{lemma}\label{lemma01} \cite{LL21} Given  $(X,\cdot)\in $ {\bf Grp}, $x\in X, A\in L^X$.

{\rm (1)} $\top_{\{x\}}\odot  A=\mathcal {L}^{\rightarrow}_x( A)=\mathcal
{L}^{\leftarrow}_{x^{-1}}( A)$.

{\rm (2)} $ A\odot \top_{\{x\}}=\mathcal
{R}^{\rightarrow}_x( A)=\mathcal {R}^{\leftarrow}_{x^{-1}}( A)$.

{\rm (3)}  $\top_{\{e\}}\odot  A= A= A\odot \top_{\{e\}}$, where $e$ is the identity element. %
%
%
%
%
%
%
%
\end{lemma}

\begin{definition} \cite{AHS} A concrete category ({\bf A}, $U$)  is called {\it topological} over {\bf X} if each $U$- source $(X\stackrel{f_i}{\longrightarrow} U(A_i))_{i\in I}$ has a unique $U$-initial lift $(A\stackrel{f_i}{\longrightarrow} A_i)_{i\in I}$, that is, an {\bf X}-morphism $f:U(B)\longrightarrow X$ is an {\bf A}-morphism iff each composition $f_i\circ f: U(B)\longrightarrow U(A_i)$ is an {\bf A}-morphism.
\end{definition}

\begin{theorem} \cite{AHS} A concrete category ({\bf A},$U$)  is  {\it topological} over {\bf X} iff each $U$-sink $(U(A_i)\stackrel{f_i}{\longrightarrow} X )_{i\in I}$ has a unique $U$-final lift $( A_i\stackrel{f_i}{\longrightarrow}A)_{i\in I}$, that is, an {\bf X}-morphism $f:X\longrightarrow U(B)$ is an {\bf A}-morphism iff each composition $f \circ f_i : U(A_i)\longrightarrow U(B)$ is an {\bf A}-morphism.
\end{theorem}

\section{The category of stratified $L$-convex spaces}

This section will show that the category of stratified $L$-convex spaces forms a topological category over {\bf SET} by demonstrating that each
$U$-source  has a unique initial  lift. Based on this, the concepts of subspaces, join spaces, quotient spaces, and product spaces of stratified $L$-convex spaces will be defined.

The family $\{ A_i\}_{i\in I}\subseteq L^X$ is referred  directed if for every $ A_{i_1},  A_{i_2}$, there exists $ A_{i_3}$ s.t.  $  A_{i_1}\leq A_{i_3}$ and $ A_{i_2}\leq A_{i_3} $. We use  $\{ A_i\}_{i\in I}\subseteq^\uparrow  L^X$ to denote a  directed subfamily of $L^X$, and use ${\bigvee\limits_{i\in
I}}^\uparrow A_i$ to denote its supremum.

\begin{definition} \label{def3} (\cite{ZL20}) An $L$-convex structure on  $X$ is a collection $\mathbb{C}\subseteq L^X$ fulfills:

(1) $\hat{\top}, \hat{\bot}\in \mathbb{C}$,

(2)  $\bigwedge\limits_{t\in T} A_t\in \mathbb{C}$ for any $\{ A_t\}_{t\in T}\subseteq \mathbb{C}$,

(3)   ${\bigvee\limits_{t\in T}}^\uparrow A_t\in \mathbb{C}$ for any $\{ A_t\}_{t\in T}\subseteq^\uparrow  \mathbb{C}$.

Then $(X, \mathbb{C})$ is termed an $L$-convex space ($L$-cvs for short). Furthermore, it is termed stratified whenever $\hat{a}\in \mathbb{C}$ for any $a\in L$.
%
\end{definition}

A mapping  $(X,\mathbb{C})\stackrel{f}\longrightarrow (Y,\mathbb{G})$ between $L$-cvss is termed $L$-convexity-preserving (abbreviated as LCP)
 whenever $ A\in \mathbb{G}$ ensures $f^\leftarrow(A)\in \mathbb{C}$.

 Obviously, the category of stratified $L$-cvss along with
LCPs, denoted as {\bf SLCVS}, forms a concrete category over {\bf SET} concerning the forgetful functor $U$.

Note that when $L=\{\bot,\top\}$, {\bf SLCVS} degenerates into {\bf CVS}.

In \cite{XZY202}, Xiu introduced the notion of $L$-convex remote-neighborhood system and shown that $L$-cvs can be described by this notion. One can redefine it in stratified $L$-cvs similarly.

\begin{definition} Given a stratified $L$-cvs $(X,\mathbb{C})$ and $x_a\in J(L^X)$. The family $$R^{\mathbb{C}}_{x_a}=\{A\in L^X|\exists C\in \mathbb{C} \ s.t.\ x_a\leq C'\leq A'\}$$ is termed the $L$-convex remote-neighborhood system at $x_a$.
\end{definition}

\begin{proposition} \label{prpocpr}  A mapping  $(X,\mathbb{C})\stackrel{f}\longrightarrow (Y,\mathbb{G})$ is an LCP iff  $f^\leftarrow (A)\in R^{\mathbb{C}}_{x_a}$ for every $x_a\in J(L^X)$ with $A\in R^{\mathbb{G}}_{f(x)_a}$.
\end{proposition}

\begin{proof} $\Longrightarrow$. Let $x_a\in J(L^X)$ and $A\in R^{\mathbb{G}}_{f(x)_a}$. Then there exists a $B\in \mathbb{G}$ s.t. $f(x)_a\leq B'\leq A'$, i.e., $x_a\leq f^\leftarrow (B)'\leq f^\leftarrow (A)'$. Because $f$ is an LCP, then $f^\leftarrow (B)\in \mathbb{C}$, thus $f^\leftarrow (A)\in R^{\mathbb{C}}_{x_a}$.

$\Longleftarrow$. Let $A\in \mathbb{G}$. Then for any $x_a\leq f^\leftarrow(A)'$, one get $f(x)_a\leq A'$, so $A\in R^{\mathbb{G}}_{f(x)_a}$, hence $f^\leftarrow (A)\in R^{\mathbb{C}}_{x_a}$ from the given condition. By the definition of  $L$-convex remote-neighborhood system, there exists a $C_{x_a}\in \mathbb{C}$ s.t. $x_a\leq C_{x_a}'\leq f^\leftarrow (A)'$. It follows that
$$f^\leftarrow(A)'=\bigvee\{x_a|x_a\leq f^\leftarrow(A)'\}\leq \bigvee\{C_{x_a}'|x_a\leq f^\leftarrow(A)'\}=\Big(\bigwedge \{C_{x_a}|x_a\leq f^\leftarrow(A)'\}\Big)'\leq f^\leftarrow(A)',$$ which means $f^\leftarrow(A)=\bigwedge \{C_{x_a}|x_a\leq f^\leftarrow(A)'\}\in \mathbb{C}$. Therefore, $f$ is an LCP.
\end{proof}

 Pang \cite{PB19} proposed the theory of base and subbase in $L$-cvs. This theory can be stated to stratified $L$-cvs with only minor modifications. Since the proofs follow a similar pattern, they are omitted.

\begin{definition}  \label{defbs} \cite{PB19} Given a stratified $L$-cvs $(X,\mathbb{C})$.

(1) $\mathbb{B}\subseteq \mathbb{C}$ is termed a base of $(X,\mathbb{C})$ whenever $\forall D\in \mathbb{C}, \exists \mathbb{B}_D{\subseteq}^\uparrow \mathbb{B}, \ {\rm s.t.} \ D={\bigvee}^\uparrow  \mathbb{B}_D$.

(2) $\mathbb{S}\subseteq \mathbb{C}$ is termed a subbase of $(X,\mathbb{C})$ whenever $\mathbb{B}_{\mathbb{S}}=\Big\{\bigwedge A_i\Big | I\neq \emptyset, \forall i\in I,  A_i\in \mathbb{S} \Big\}$
forms a base.
\end{definition}

\begin{proposition} \cite{PB19} A family $\mathbb{B}\subseteq L^X$ forms a base of one  stratified $L$-cvs  iff $\mathbb{B}$ fulfills:

(LCB1) $\hat{\bot} \in \mathbb{B}$ and $\forall a\in L, a\neq \bot$, $\hat{a}=\bigvee^\uparrow  \mathbb{B}_{a}$ for some $\mathbb{B}_{a}\subseteq ^\uparrow  \mathbb{B}$.

(LCB2) For any nonempty $\{B_i|i\in I\}\subseteq \mathbb{B}$,  $\exists \{C_j|j\in J\}\subseteq^\uparrow  \mathbb{B}$ s.t. $\bigwedge\limits_{i\in I} B_i={\bigvee\limits_{j\in J}}^\uparrow  C_j$.

(LCB3) For any $\{B_j|j\in J\}\subseteq^\uparrow  L^X$, where each $B_j={\bigvee\limits_{i\in I_j}}^\uparrow C_{ij}$ and $\{C_{ij}|i\in I_j\}\subseteq^\uparrow  \mathbb{B}$,  there is $\{G_k|k\in K\}\subseteq ^\uparrow  \mathbb{B}$ s.t. ${\bigvee\limits_{k\in K}}^\uparrow  G_k={\bigvee\limits_{j\in J}}^\uparrow {\bigvee \limits_{i\in I_j}}^\uparrow C_{ij}$.
\end{proposition}

\begin{proposition} \cite{PB19} A family $\mathbb{S}\subseteq L^X$ is a subbase of one stratified $L$-cvs  iff $\mathbb{S}$ fulfills:

(LCSB1) There is $\{B_i|i\in I\}\subseteq \mathbb{S}$ s.t. $\hat{\bot}=\bigwedge\limits_{i\in I} B_i$.

(LCSB2) For any $a\in L$, $a\neq \bot$, there is $\{B_{ij}|j\in J_i\}\subseteq \mathbb{S}$ for each $i\in I$ s.t. $\{\bigwedge\limits_{j\in J_i}B_{ij}|i\in I\}\subseteq ^\uparrow  L^X$ and $\hat{a}={\bigvee\limits_{i\in I}}^\uparrow \bigwedge\limits_{j\in J_i} B_{ij}$.
\end{proposition}

\begin{definition} Let $\{\mathbb{C}_i\}_{i\in I}$ be a family of stratified $L$-convex structures on $X$. Then  the family $\bigcup\limits_{i\in I}\mathbb{C}_i$ forms a subbase of a stratified $L$-convex structure on $X$, denoted by $\bigsqcup\limits_{i\in I} \mathbb{C}_i$, called the join structure of $\{\mathbb{C}_i\}_{i\in I}$. \end{definition}

\begin{proposition} \label{prpocp} \cite{PB19} Suppose $\mathbb{B}$ (resp., $\mathbb{S}$) is a base (resp., subbase) of  $(Y,\mathbb{G})$. Then $(X,\mathbb{C})\stackrel{f}\longrightarrow (Y,\mathbb{G})$ is an LCP iff $f^\leftarrow (A)\in \mathbb{C}$ for each $A\in \mathbb{B}$ (resp., $A\in \mathbb{S}$).
\end{proposition}

\begin{definition} \label{defcvin}  A $U$-source $(X\stackrel{f_j}{\longrightarrow} U(X_j,\mathbb{C}_j))_{j\in J}$ in {\bf SLCVS} is a collection of mappings $X\stackrel{f_j}{\longrightarrow} X_j$ with stratified $L$-cvss $(X_j,\mathbb{C}_j)$. The given source is termed to possess a unique $U$-initial lift provided there exists a stratified $L$-convex structure $\mathbb{C}$ on $X$ fulfills:

(1) Each $f_j:(X,\mathbb{C})\longrightarrow (X_j,\mathbb{C}_j)$ ($j\in J$) is an LCP,

(2)  $\forall (Y,\mathbb{G})\in ${\bf SLCVS},  $(Y,\mathbb{G})\stackrel{ f}\longrightarrow(X,\mathbb{C})$ is an LCP iff all $(Y,\mathbb{G})\stackrel{f_j\circ f}\longrightarrow (X_j,\mathbb{C}_j)$ are LCPs.
\end{definition}

\begin{proposition}  \label{tinitial} Any $U$-source $(X\stackrel{f_k}{\longrightarrow} U(X_k,\mathbb{C}_k))_{k\in K}$ in {\bf SLCVS} possesses  a unique $U$-initial lift. Hence {\bf SLCVS} is a topological category over {\bf SET}.
\end{proposition}

\begin{proof} Take $k\in K$, the family $f_k^\leftarrow (\mathbb{C}_k)=\{f_k^\leftarrow (C_k)|C_k\in \mathbb{C}_k\}$ forms a stratified $L$-convex structure on $X$. Let $\mathbb{C}=\bigsqcup \{f_k^\leftarrow (\mathbb{C}_k)|k\in K\}$, i.e., the stratified $L$-convex structure generated by  the family $\bigcup \{f^{\leftarrow}_k(\mathbb{C}_k)|k\in K\}$ as a subbase. It concludes directly from Definition \ref{defbs} and Proposition \ref{prpocp} that $(X,\mathbb{C})$ fulfills the $U$-initial lift condition.
\end{proof}

\begin{definition} The stratified $L$-convex structure $\mathbb{C}=\bigsqcup \{f_j^\leftarrow (\mathbb{C}_j)|j\in J\}$ defined in Proposition \ref{tinitial} is called the initial structure of the $U$-source $(X\stackrel{f_j}{\longrightarrow} U(X_j,\mathbb{C}_j))_{j\in J}$.
\end{definition}

As is well known, the substructure and product structure are special types of initial structures. Therefore, Proposition \ref{tinitial} implies the following conclusions.

\begin{corollary}  Let $(X,\mathbb{C})$ belong to {\bf SLCVS}, with $Y\subseteq X$ and $i:Y\longrightarrow X$ is the inclusion mapping. Then the $U$-source $Y\stackrel{i}\longrightarrow U(X,\mathbb{C})$ has an initial structure, denoted by $\mathbb{C}|_Y$. The pair $(Y,\mathbb{C}|_Y)$ is called a subspace of  $(X,\mathbb{C})$.
\end{corollary}

\begin{corollary}  Assume that $(X_i,\mathbb{C}_i)$ ($i\in I$) belong to {\bf SLCVS} and let $p_i:\prod\limits_{i\in I} X_i\longrightarrow X_i$ denote the projection mappings. Then the $U$-source $\big(\prod\limits_{i\in I} X_i\stackrel{p_i}\longrightarrow U(X_i,\mathbb{C}_i)\big)_{i\in I}$ possesses an initial structure, denoted by $\prod\limits_{i\in I}\mathbb{C}_i$. The pair $(\prod\limits_{i\in I} X_i,\prod\limits_{i\in I}\mathbb{C}_i)$  is called the product space of $(X_i,\mathbb{C}_i)$.
\end{corollary}

Pang \cite{PB19} provided a concise characterization of the finite product.

\begin{proposition} \label{propfop} Let $(X_i,\mathbb{C}_i)$ ($i=1,2, \ldots, n$) belong to {\bf SLCVS}. Then  the product  $\prod\limits_{i=1}^n\mathbb{C}_i=\{C_1\times C_2\times \ldots \times C_n|C_i\in \mathbb{C}_i\}$.
\end{proposition}

\begin{definition} A $U$-sink $(U(X_j,\mathbb{C}_j)\stackrel{f_j}{\longrightarrow} X)_{j\in J}$ in {\bf SLCVS} is a collection of mappings $X_j\stackrel{f_j}{\longrightarrow} X$ with stratified $L$-cvss $(X_j,\mathbb{C}_j)$. The given sink is termed to possess a unique $U$-final lift provided there exists a stratified $L$-convex structure $\mathbb{C}$ on $X$ fulfills:

(1) Each $f_j:(X_j,\mathbb{C}_j)\longrightarrow(X,\mathbb{C})$ ($j\in J$) is an LCP,

(2)  $\forall (Y,\mathbb{G})\in ${\bf SLCVS},  $(X,\mathbb{C})\stackrel{f}\longrightarrow(Y,\mathbb{G})$ is  LCP iff all $(X_j,\mathbb{C}_j)\stackrel{f\circ f_j}\longrightarrow (Y,\mathbb{G})$ are LCPs.
\end{definition}
%

\begin{proposition} \cite{PB232} \label{tfinal} Each $U$-sink $(U(X_j,\mathbb{C}_j)\stackrel{f_j}{\longrightarrow} X)_{j\in J}$ in {\bf SLCVS} possesses a unique $U$-final lift defined by $\mathbb{C}=\{C\in L^X|\forall j\in J, f_j^\leftarrow (C)\in \mathbb{C}_j\}$.
\end{proposition}


\begin{definition} The stratified $L$-convex structure $\mathbb{C}$ defined in Proposition \ref{tfinal} is referred to as the final structure of the given sink.
\end{definition}

As is well known, the quotient structures are special types of initial structures. Therefore, Proposition \ref{tfinal} implies the following conclusion.

\begin{corollary}  \label{coroqu} Let $(X,\mathbb{C})$ belong to {\bf SLCVS}, with $f:X\longrightarrow Y$ being a surjective mapping, and let $\mathbb{C}_f$ denote the final structure of the $U$-sink   $f:U(X,\mathbb{C})\longrightarrow Y$. Then, the pair $(Y,\mathbb{C}_f)$ is called the quotient space of  $(X,\mathbb{C})$ w.r.t $f$.
\end{corollary}

\section{The category of stratified $L$-convex groups}

This section investigates a category of stratified $L$-convex groups, along with some of their properties and characterizations.

Given a group $(X,\cdot)$, we use $m: X\times X\longrightarrow X$ to represent  the multiplication, and $r:X\longrightarrow X$ to the inverse operation.

\begin{definition} Given a stratified  $L$-convex structure $\mathbb{C}$ on a group $(X,\cdot)$, the triple
$(X,\cdot,\mathbb{C})$ is termed a stratified  $L$-convex
group provided that  both the  multiplication  $(X\times X, \mathbb{C}\times
\mathbb{C})\stackrel{m}{\longrightarrow} (X,\mathbb{C})$ and the inverse operation $(X,\mathbb{C})\stackrel{r}{\longrightarrow}
(X,\mathbb{C})$  are LCPs. \end{definition}

When $L=\{\bot,\top\}$, stratified $L$-convex spaces reduce to  convex spaces, and  LCPs reduce to CPs. This motivates us to define the notion below.

\begin{definition} Given a convex structure $\mathcal{C}$ on a group $(X,\cdot)$,  the triple $(X,\cdot,\mathcal{C})$ is termed a convex
group provided that  both the multiplication $m$ and  the inverse operation $r$ are CPs.
\end{definition}

Let {\bf SLCG} (resp., {\bf CG}) denote the category of stratified $L$-convex groups (resp., convex groups) and LCP (resp., CP) group homomorphisms.

Obviously, {\bf SLCG} forms a concrete category over {\bf Grp} concerning the forgetful functor $U$.

\begin{example} (1) Given a group $(X,\cdot)$, then both $\mathbb{C}_{0}=\{\hat{a}|a\in L\}$  and  $\mathbb{C}_1=L^X$,  together with $(X,\cdot)$,  forms a stratified $L$-convex group, respectively.

(2) Any convex group generates a stratified $L$-convex group, as stated in Proposition \ref{thmgg} below.
\end{example}

Next, we discuss the properties and characterizations of {\bf SLCG}. As a corollary, {\bf CG} also leads to similar conclusions, which we will not restate separately.

\begin{lemma} \label{lemmacv} Given a stratified
$L$-cvs $(X,\mathbb{C})$ and $x\in X$, we define
$\varphi_x:X\longrightarrow X\times X$ by $u\mapsto (x,u)$. Then $\varphi_x$ is an LCP.
\end{lemma}

\begin{proof} Obviously, $p_2\circ \varphi_x=id_X$ is an LCP. Additionally, $p_1\circ \varphi_x$ is an LCP because it always takes the element $x$. Therefore, $\varphi_x$ is an LCP through the initial property of product structure.
\end{proof}

\begin{remark} Lemma \ref{lemmacv} does not hold for $L$-cvs since the constant mapping between non-stratified $L$-cvs is not necessarily an LCP in general.
\end{remark}

\begin{proposition} Let $(X,\cdot,\mathbb{C})$ be a stratified
$L$-cvs and $x\in X$. Then the left (resp., right) translation $\mathcal {L}_x$ (resp., $\mathcal {R}_x$) is an LCP homeomorphism.
\end{proposition}

\begin{proof} From $\mathcal {L}_x=m\circ \varphi_x$ and $\mathcal {L}_{x^{-1}}=m\circ \varphi_{x^{-1}}$ we know that $\mathcal {L}_x$ and $\mathcal {L}_{x^{-1}}$ are LCPs. It
follows that $\mathcal {L}_x$ is a homeomorphism by $\mathcal
{L}_x\circ \mathcal {L}_{x^{-1}}=\mathcal {L}_{x^{-1}}\circ \mathcal
{L}_x=id_X$. The case of $\mathcal {R}_x$ is similar.
\end{proof}

\begin{theorem} (Localizable Theorem) \label{theorem01} For an $L$-convex group $(X,\cdot,\mathbb{C})$, $A\in R^\mathbb{C}_{x_a}\Longleftrightarrow
\top_{\{x^{-1}\}}\odot A\in R^\mathbb{C}_{e_a}\Longleftrightarrow A\odot \top_{\{x^{-1}\}}\in
R^\mathbb{C}_{e_a}$.
\end{theorem}
\begin{proof} We only demonstrate $A\in R^\mathbb{C}_{x_a}\Longleftrightarrow
\top_{\{x^{-1}\}}\odot A\in R^\mathbb{C}_{e_a}$, and the rest follows in a similar manner. From Lemma \ref{lemma01}, it concludes that $\mathcal
{L}_{x^{-1}}^\leftarrow(A)=\top_{\{x\}}\odot A, \mathcal
{L}^\leftarrow_x(A)=\top_{\{x^{-1}\}}\odot A$. Since
$\mathcal {L}_x$ and $\mathcal {L}_{x^{-1}}$ are LCPs, we can deduce the following results from Proposition \ref{prpocpr}
$$A\in  R^\mathbb{C}_{x_a}= R^\mathbb{C}_{\mathcal {L}_{x}(e)_a}\Longrightarrow \mathcal
{L}_{x}^{\leftarrow}(A)=\top_{\{x^{-1}\}}\odot A\in R^\mathbb{C}_{e_a},$$
$$\top_{\{x^{-1}\}}\odot A\in R^\mathbb{C}_{e_a}\Longrightarrow \top_{\{x^{-1}\}}\odot A\in R^\mathbb{C}_{(\mathcal {L}_{x^{-1}}(x))_a}\Longrightarrow
A=\mathcal {L}_{x^{-1}}^{\leftarrow}(\top_{\{x^{-1}\}}\odot A)\in
R^\mathbb{C}_{x_a}.\qedhere$$
\end{proof}

The localization theorem states that an $L$-convex group is completely determined by its remote neighborhood systems at the identity element $e$.

For two mappings $X_i\stackrel{f_i}\longrightarrow Y_i$ $(i=1,2)$, we define one mapping $X_1\times X_2\stackrel{f_1\times f_2}\longrightarrow Y_1\times Y_2$ by $(x_1,x_2)\mapsto (f_1(x_1),f_2(x_2))$.

\begin{lemma} \label{lempr} Let $(X_i,\mathbb{C}_i)\stackrel{f_i}\longrightarrow (Y_i,\mathbb{G}_i)$ $(i=1,2)$ be LCPs between stratified $L$-cvss. Then  $f_1\times f_2$ is also an LCP between the product spaces of stratified $L$-cvss.
\end{lemma}
\begin{proof} From Proposition \ref{propfop} one conclude that the products $\mathbb{C}_1\times \mathbb{C}_2=\{C_1\times C_2|C_i\in \mathbb{C}_i\}$ and $\mathbb{G}_1\times \mathbb{G}_2=\{G_1\times G_2|G_i\in \mathbb{G}_i\}$. For any $G_i\in \mathbb{G}_i$, one obtain that $f_i^\leftarrow(G_i)\in \mathbb{C}_i$ by $f_i$ is an LCP. It follows that $(f_1\times f_2)^\leftarrow(G_1\times G_2)=f_1^\leftarrow(G_1)\times f_2^\leftarrow(G_2) \in \mathbb{C}_1\times \mathbb{C}_2$. Hence $f_1\times f_2$ is an LCP.
\end{proof}

For a group $(X,\cdot)$, we define $k:X\times X\longrightarrow X$  by $(x,y)\mapsto xy^{-1}$.

\begin{lemma} \label{lempro} Let $(X,\cdot)$ be a group and $A,B\in L^X$. Then $k^\rightarrow (A\times B)=A\odot B^{-1}$.
\end{lemma}

\begin{proof} For any $A,B\in L^X$ and $z\in X$, \begin{eqnarray*}k^\rightarrow (A\times B)(z)&=&\bigvee_{xy^{-1}=z}\Big(A(x)\wedge B(y)\Big)=\bigvee_{xy^{-1}=z}\Big(A(x)\wedge B^{-1}(y^{-1})\Big)\\
&=&\bigvee_{xw=z}\Big(A(x)\wedge B^{-1}(w)\Big)=(A\odot B^{-1})(z),\end{eqnarray*}which means $k^\rightarrow (A\times B)=A\odot B^{-1}$.
\end{proof}

\begin{proposition} \label{theoremch} Give a group $(X,\cdot)$ and a stratified $L$-cvs $(X,\mathbb{C})$, $(X,\cdot,\mathbb{C})\in $ {\bf SLCG}  iff the mapping $k$ is an LCP.
\end{proposition}

\begin{proof} If $(X,\cdot,\mathbb{C})\in ${\bf SLCG}, then both $m$ and $r$ are LCPs. Since $id_{X}$ is an LCP and by Lemma \ref{lempr}, we conclude that $id_{X}\times r$ is an LCP. Hence, $k=m\circ (id_{X}\times r)$ is also an LCP.

Conversely, assume that $k$ is an LCP. From Lemma \ref{lemmacv} we know that $\varphi_e$ is an LCP, so $r=k\circ \varphi_e$ is an LCP. By Lemma \ref{lempr}, it holds that $id_{X}\times r$ is an LCP. Hence,  $m=k\circ (id_{X}\times r)$ is an LCP.\end{proof}

Similar to topological groups, the next theorem demonstrates that stratified $L$-convex groups can be characterized by the operation $\odot$.
\begin{theorem} \label{thmred} The triple $(X,\cdot,\mathbb{C})\in $ {\bf SLCG}  iff for any $x,y\in X$ and $A\in R^{\mathbb{C}}_{{xy^{-1}}_a}$, there exists $B,C\in \mathbb{C}$ with $B\in R^{\mathbb{C}}_{x_a}$ or $C\in R^{\mathbb{C}}_{y_a}$ s.t. $(B'\odot \hat{\top}^{-1})\vee (\hat{\top}\odot C'^{-1}) \leq A'$.
\end{theorem}

\begin{proof} $\Longrightarrow$. Let $x,y\in X$, $A\in R^{\mathbb{C}}_{{xy^{-1}}_a}=R^{\mathbb{C}}_{k(x,y)_a}$. Since $k$ is an LCP, we have  $k^{\leftarrow}(A)\in R^{\mathbb{C}\times \mathbb{C}}_{{(x,y)}_a}$. By Proposition \ref{propfop}, there exists $B,C\in \mathbb{C}$ s.t. $(x,y)_a\leq (B\times C)'\leq k^{\leftarrow}(A)'$. It is easily seen that $(B\times C)'=(B'\times \hat{\top}) \vee (\hat{\top}\times C')$, so $$(x,y)_a\leq \Big((B'\times \hat{\top}) \vee (\hat{\top}\times C')\Big)\leq k^{\leftarrow}(A)'.$$ By $a\in J(L)$, one have $$(x,y)_a\leq (B'\times \hat{\top}) \ {\rm or}\  (x,y)_a\leq (\hat{\top}\times C')\Longrightarrow x_a\leq B' \ {\rm or}\ y_a\leq C'\Longrightarrow B\in R^{\mathbb{C}}_{x_a}  \ {\rm or}\  C\in R^{\mathbb{C}}_{y_a}.$$Furthermore,  \begin{eqnarray*}(B'\times \hat{\top}) \vee (\hat{\top}\times C')\leq k^{\leftarrow}(A')&\Longleftrightarrow& k^\rightarrow\Big((B'\times \hat{\top}) \vee (\hat{\top}\times C')\Big)\leq A'\\
&\Longleftrightarrow& k^\rightarrow(B'\times \hat{\top}) \vee k^\rightarrow(\hat{\top}\times C')\leq A'\\
&\stackrel{\rm Lemma \ref{lempro}}{\Longleftrightarrow}& (B'\odot \hat{\top}^{-1}) \vee (\hat{\top}\odot C'^{-1})\leq A',\end{eqnarray*}as desired.

$\Longleftarrow$. Assume the given condition is satisfied. Let $x,y\in X$, $A\in R^{\mathbb{C}}_{{xy^{-1}}_a}=R^{\mathbb{C}}_{k(x,y)_a}$. Then there exists $B,C\in \mathbb{C}$ with $B\in R^{\mathbb{C}}_{x_a}$ or $C\in R^{\mathbb{C}}_{y_a}$ s.t. $(B'\odot \hat{\top}^{-1})\vee (\hat{\top}\odot C'^{-1}) \leq A'$.

From $B,C\in \mathbb{C}$ and $B\in R^{\mathbb{C}}_{x_a}$ or $C\in R^{\mathbb{C}}_{y_a}$, one get $x_a\leq B'$ or $y_a\leq  C'$. This means $$(x,y)_a\leq \Big((B'\times \hat{\top}) \vee (\hat{\top}\times C')\Big)=(B\times C)'.$$

From $(B'\odot \hat{\top}^{-1})\vee (\hat{\top}\odot C'^{-1}) \leq A'$, one get $$k^\rightarrow(B'\times \hat{\top}) \vee k^\rightarrow(\hat{\top}\times C')\leq A'\Longrightarrow (B'\times \hat{\top}) \vee (\hat{\top}\times C')\leq k^{\leftarrow}(A)'.$$ Together above, we obtain $(x,y)_a\leq (B\times C)'\leq k^{\leftarrow}(A)'$, i.e., $k^{\leftarrow}(A)\in R^{\mathbb{C}\times \mathbb{C}}_{(x,y)_a}$. It follows by Proposition \ref{prpocpr}  that $k$ is an LCP, as desired.\end{proof}

\begin{remark} \label{remrn} Theorems \ref{theorem01} and \ref{thmred} differ from the corresponding results in topological groups, the former employs ``remote neighborhoods'' while the latter uses ``neighborhoods'' as the tools of investigation.
\end{remark}

\section{ The initial and final structure in {\bf SLCG}}

This section discusses the initial and final structure in {\bf SLCG}, and shows that {\bf SLCG} is a topological category over {\bf Grp}.

\begin{definition} \label{defgin} A $U$-source $((X,\cdot)\stackrel{f_j}{\longrightarrow} U(X_j,\cdot,\mathbb{C}_j))_{j\in J}$ in {\bf SLCG} over {\bf Grp} is a collection of group homomorphisms $(X,\cdot)\stackrel{f_j}{\longrightarrow} (X_j,\cdot)$ with stratified $L$-convex groups $(X_j,\cdot,\mathbb{C}_j)$. The given source is termed to possess a unique $U$-initial lift provided there is a stratified $L$-convex group $(X,\cdot,\mathbb{C})$ fulfills:

(1) Each $f_j:(X,\cdot,\mathbb{C})\longrightarrow (X_j,\cdot,\mathbb{C}_j)$ ($j\in J$) is an LCP group homomorphism,

(2)  $\forall (Y,\cdot,\mathbb{G})\in ${\bf SLCG}, a group homomorphism $(Y,\cdot,\mathbb{G})\stackrel{ f}\longrightarrow(X,\cdot,\mathbb{C})$ is an LCP iff all $(Y,\cdot,\mathbb{G})\stackrel{f_j\circ f}\longrightarrow (X_j,\cdot,\mathbb{C}_j)$ are LCPs.
\end{definition}

\begin{theorem} \label{tngin} Each $U$-source $\big((X,\cdot) \stackrel{f_j}{\longrightarrow} (X_j,\cdot,\mathbb{C}_j)\big)_{j\in J}$ of {\bf SLCG} possesses  a unique $U$-initial lift. Hence {\bf SLCG} is a topological category over {\bf Grp}.
\end{theorem}

\begin{proof} Let $\big((X,\cdot) \stackrel{f_j}{\longrightarrow} (X_j,\cdot,\mathbb{C}_j)\big)_{j\in J}$ be a $U$-source in {\bf SLCG}. By  Proposition \ref{tinitial},  $\big(X \stackrel{f_j}{\longrightarrow} (X_j,\mathbb{C}_j)\big)_{j\in J}$ forms a $U$-source  in {\bf SLCVS}, with the initial structure $\mathbb{C}=\bigsqcup \{f_j^\leftarrow (\mathbb{C}_j)| j\in J\}$. We verify below that $(X,\cdot, \mathbb{C})$ is the $U$-initial lift of the source $\big((X,\cdot) \stackrel{f_j}{\longrightarrow} (X_j,\cdot_j,\mathbb{C}_j)\big)_{j\in J}$ in {\bf SLCG}.

Since $\mathbb{C}$ fulfills the $U$-initial condition in Definition \ref{defcvin}, $(X,\cdot, \mathbb{C})$ also fulfills that condition in Definition \ref{defgin}. Hence, it suffices to verify that $(X,\cdot, \mathbb{C})$ constitutes a stratified $L$-convex group.

It is easy to see that the following diagram is commutative.
$$\bfig \morphism(0,0)|a|/@{<-}@<0pt>/<800,0>[X_j`
X;f_j] \morphism(0,0)|l|/@{<-}@<0pt>/<0,-400>[X_j`X_j\times X_j;k_j]
\morphism(800,0)|r|/@{<-}@<0pt>/<0,-400>[X`X\times X;k] \morphism(0,-400)|a|/@{<-}@<0pt>/<800,0>[X_j\times X_j`X\times X;f_j\times f_j] \efig
$$
 From Lemma \ref{lempr} and Proposition \ref{theoremch}, we know that for any $j\in J$,  both $f_j\times f_j$ and $k_j$ are LCPs. Therefore, $f_j\circ k=k_j\circ (f_j\times f_j)$ is also an LCP. Hence  $k$ is an LCP, since $ \mathbb{C}$ fulfills the $U$-initial lift condition (2) in Definition \ref{defcvin}. Therefore, according to Proposition \ref{theoremch}, $(X,\cdot, \mathbb{C})$ is a stratified $L$-convex group. \end{proof}




%
%

\begin{definition} For a $U$-source $((X,\cdot)\stackrel{f_j}{\longrightarrow} U(X_j,\cdot,\mathbb{C}_j))_{j\in J}$ in {\bf SLCG}, the stratified $L$-convex group $(X,\cdot,\mathbb{C})$ defined in Theorem \ref{tngin} is termed the initial structure of the given source.
\end{definition}

Theorem \ref{tngin} implies the following conclusions.

\begin{corollary}  Let $(X,\cdot,\mathbb{C})$ belong to {\bf SLCG}, and $Y$ be a subgroup of $X$. Then the triple $(Y,\cdot, \mathbb{C}|_Y)$ is a stratified $L$-convex group, called a subspace of $(X,\cdot,\mathbb{C})$.
\end{corollary}

\begin{corollary}  Let $(X_i,\cdot,\mathbb{C}_i)$ ($i\in I$) belong to {\bf SLCG}. Then  $(\prod\limits_{i\in I} X_i,\cdot, \prod\limits_{i\in I}\mathbb{C}_i)$ belongs to {\bf SLCG}, called the product space of $(X_i,\cdot,\mathbb{C}_i)$.
\end{corollary}

\begin{corollary}  Let $(X,\cdot,\mathbb{C}_i)$ ($i\in I$) belong to {\bf SLCG}. Then  $(X,\cdot, \bigsqcup\mathbb{C}_i)$ belongs to {\bf SLCG}, called the join space of $(X,\cdot,\mathbb{C}_i)$.
\end{corollary}

\begin{definition} A $U$-sink $(U(X_j,\cdot,\mathbb{C}_j)\stackrel{f_j}{\longrightarrow} (X,\cdot))_{j\in J}$ in {\bf SLCG} over {\bf Grp} is a collection of group homomorphisms $(X_j,\cdot)\stackrel{f_j}{\longrightarrow} (X,\cdot)$ with stratified $L$-convex groups $(X_j,\cdot,\mathbb{C}_j)$. The given sink is termed to possess a unique  $U$-initial lift provided there is a stratified $L$-convex group $(X,\cdot,\mathbb{C})$ fulfills:

(1) Each $f_j:(X_j,\cdot,\mathbb{C}_j)\longrightarrow(X,\cdot,\mathbb{C})$ ($j\in J$) is an LCP group homomorphism,

(2)  $\forall (Y,\cdot,\mathbb{G})\in ${\bf SLCG}, a group homomorphism $(X,\cdot,\mathbb{C})\stackrel{ f}\longrightarrow(Y,\cdot,\mathbb{G})$ is an LCP iff all $(X_j,\cdot,\mathbb{C}_j)\stackrel{f\circ f_j}\longrightarrow(Y,\cdot,\mathbb{G})$ are LCPs.\end{definition}


\begin{proposition} \label{prpocp2}  $(X,\mathbb{C})\stackrel{f}\longrightarrow (Y,\mathbb{G})$ is an LCP iff $\forall A'\in \mathbb{G}$, $[f^\leftarrow (A)]'\in \mathbb{C}$.
\end{proposition}

\begin{proof} The necessity is clear, only sufficiency remains for proving. Indeed, let $A\in \mathbb{G}$, i.e., $A=(A')'\in \mathbb{G}$, then from the given condition, one obtain $[f^\leftarrow(A')]'\in \mathbb{C}$, i.e.,  $f^\leftarrow(A)\in \mathbb{C}$, as desired.
\end{proof}

\begin{definition} Let $(X,\mathbb{C})\stackrel{f}\longrightarrow (Y,\mathbb{G})$ be a mapping between $L$-cvss. Then $f$ is called complement $L$-convex-to-convex if $[f^\rightarrow (A)]'\in \mathbb{G}$ for all $A'\in \mathbb{C}$.
\end{definition}

\begin{theorem} \label{tngfinal} Each $U$-sink $\big((X_j,\cdot,\mathbb{C}_j) \stackrel{f_j}{\longrightarrow}(X,\cdot)\big)_{j\in J}$ in {\bf SLCG} possesses  a unique $U$-final lift, provided that each $f_j$ is a surjective and complement $L$-convex-to-convex mapping.
\end{theorem}

\begin{proof} Let  $\big((X_j,\cdot,\mathbb{C}_j) \stackrel{f_j}{\longrightarrow}(X,\cdot)\big)_{j\in J}$ be a $U$-sink in {\bf SLCG}. By Proposition \ref{tfinal}, $\big((X_j,\mathbb{C}_j) \stackrel{f_j}{\longrightarrow}X\big)_{j\in J}$ forms a sink  in {\bf SLCVS}, with the final structure  $\mathbb{C}=\{C\in L^X|\forall j\in J, f_j^\leftarrow (C)\in \mathbb{C}_j\}$. We examine  that $(X,\cdot, \mathbb{C})$ is the final lift of $\big((X_j,\cdot,\mathbb{C}_j) \stackrel{f_j}{\longrightarrow}(X,\cdot)\big)_{j\in J}$.

We only need to verify that $(X,\cdot, \mathbb{C})$ forms an  $L$-convex group by showing that $k$ is an LCP. According to Proposition \ref{prpocp2},  we require to examine $\forall A'\in \mathbb{C}$, $[k^\leftarrow (A)]'\in \mathbb{C}\times \mathbb{C}$.



It is easily observed that for any $j\in J$,  the following diagram is commutative.
$$\bfig \morphism(0,0)|a|/@{->}@<0pt>/<800,0>[X_j`
X;f_j] \morphism(0,0)|l|/@{<-}@<0pt>/<0,-400>[X_j`X_j\times X_j;k_j]
\morphism(800,0)|r|/@{<-}@<0pt>/<0,-400>[X`X\times X;k] \morphism(0,-400)|a|/@{->}@<0pt>/<800,0>[X_j\times X_j`X\times X;f_j\times f_j] \efig
$$Since $f_j$ and $k_j$ are LCPs, we have $k\circ (f_j\times f_j)=f_j\circ k_j$ is LCP. It follows that for any $A'\in \mathbb{C}$,
$$(f_j\times f_j)^\leftarrow\circ k^\leftarrow(A')=(f_j\circ k_j)^\leftarrow (A')\in \mathbb{C}_j\times \mathbb{C}_j.$$ According to Proposition \ref{propfop}, there exist $B_j, G_j\in \mathbb{C}_j$ s.t. $$(f_j\times f_j)^\leftarrow\circ k^\leftarrow(A')=(f_j\circ k_j)^\leftarrow (A')=B_j\times G_j,$$ then from $f_j$ is surjective mapping, we get \begin{eqnarray*}& &k^\leftarrow(A)=(f_j\times f_j)^\rightarrow([B_j\times G_j]')\\
&\Longrightarrow& k^\leftarrow(A)=(f_j\times f_j)^\rightarrow\Big((B_j'\times \hat{\top}) \vee (\hat{\top}\times G_j')\Big)\\
&\Longrightarrow&k^\leftarrow(A)=(f_j\times f_j)^\rightarrow(B_j'\times \hat{\top}) \vee (f_j\times f_j)^\rightarrow(\hat{\top}\times G_j')\\
&\Longrightarrow&k^\leftarrow(A)=(f_j^\rightarrow(B_j')\times f_j^\rightarrow(\hat{\top})) \vee (f_j^\rightarrow(\hat{\top})\times f_j^\rightarrow(G_j'))\\
&\Longrightarrow&k^\leftarrow(A)=(f_j^\rightarrow(B_j')\times \hat{\top}) \vee (\hat{\top}\times f_j^\rightarrow(G_j'))\\
&\Longrightarrow&k^\leftarrow(A)=[f_j^\rightarrow(B_j')'\times  f_j^\rightarrow(G_j')']'\\
&\Longrightarrow&[k^\leftarrow(A)]'=k^\leftarrow(A')=f_j^\rightarrow(B_j')'\times  f_j^\rightarrow(G_j')'.\end{eqnarray*}From $B_j=(B_j')', G_j=(G_j')'\in \mathbb{C}_j$ and that $f_j$ is complement $L$-convex-to-convex mapping, one have $f_j^\rightarrow(B_j')', f_j^\rightarrow(G_j')'\in \mathbb{C}$, so $[k^\leftarrow(A)]'=f_j^\rightarrow(B_j')'\times  f_j^\rightarrow(G_j')'\in \mathbb{C}\times \mathbb{C}$.
\end{proof}

\begin{proposition} \label{propqu} Let $(X,\cdot, \mathbb{C})$ belong to {\bf SLCG}, $N\subseteq X$ a normal subgroup of $X$, $X^*=X/N=\{zN|z\in X\}$ and  $q:X\longrightarrow X^*, z|\longrightarrow zN$ be the natural mapping. Then $q$ is a surjective complement $L$-convex-to-convex mapping between $(X, \mathbb{C})$ and the quotient space $(X^*, \mathbb{C}_q)$.
\end{proposition}

\begin{proof} Obviously, $q$ is a surjective mapping. We examine below that $q$ is an complement $L$-convex-to-convex mapping.

 Let $A'\in \mathbb{C}$. Then for any $x\in X$
 \begin{eqnarray*}q^\leftarrow (q^\rightarrow(A))(x)&=&q^\rightarrow(A)(xN)=\bigvee_{yN=xN}A(y)\\
 &=&\bigvee_{x^{-1}y\in N}A(y)=\bigvee_{z=x^{-1}y\in N}A(y)\\
 &=&\bigvee_{z\in N}A(xz)=\bigvee_{z\in N}\mathcal {R}^{\leftarrow}_{z}(A)(x)\\
 &=&\Big(\bigvee_{z\in N}\mathcal {R}^{\leftarrow}_{z}(A)\Big)(x).
 \end{eqnarray*}This means that $\Big(q^\leftarrow (q^\rightarrow(A))\Big)'=\Big(\bigvee\limits_{z\in N}\mathcal {R}^{\leftarrow}_{z}(A)\Big)'$, so $q^\leftarrow ([q^\rightarrow(A)]')=\bigwedge\limits_{z\in N}\mathcal {R}^{\leftarrow}_{z}(A')$. Because each $\mathcal {R}_{z}$ is an LCP and $A'\in \mathbb{C}$, we get $\mathcal {R}^{\leftarrow}_{z}(A')\in \mathbb{C}$, then $q^\leftarrow ([q^\rightarrow(A)]')=\bigwedge\limits_{z\in N}\mathcal {R}^{\leftarrow}_{z}(A')\in \mathbb{C}$. It follows by the definition of quotient space $(X^\ast, \mathbb{C}_q)$ that $[q^\rightarrow(A)]'\in \mathbb{C}_q$. Hence $q$ is an complement $L$-convex-to-convex mapping.
\end{proof}

\begin{corollary} Let $(X,\cdot, \mathbb{C})$ belong to {\bf SLCG}, $N\subseteq X$ a normal subgroup of $X$, $X^*=X/N=\{zN|z\in X\}$ and  $q:X\longrightarrow X^*$ be the natural mapping. Then $(X^*, \cdot, \mathbb{C}_q)$ belong to {\bf SLCG}, called a quotient space of $(X,\cdot, \mathbb{C})$.
\end{corollary}

\begin{proof} It follows by Corollary \ref{coroqu}, Theorem \ref{tngfinal} and Proposition \ref{propqu}.
\end{proof}

\begin{remark} Clearly, {\bf CG} also has initial and final structures, as well as subspaces, join spaces, product spaces, and quotient spaces. We will continue using the notation and terminology of {\bf SLCG} and will not introduce new symbols or terms to state these results.
\end{remark}

\begin{remark} \label{remfinal} Theorem \ref{tngfinal} and Proposition \ref{propqu} differ from the corresponding results in topological groups. The former assumes the mappings to be complement $L$-convex-to-convex, while the latter assumes them to be open mappings. Furthermore, in the study of convex spaces, complement convex-to-convex mappings are not considered; instead, convex-to-convex mappings are typically used.
\end{remark}

\section{The categorical relationships between {\bf CG} and {\bf SLCG}.}

This section will shows that {\bf CG} can
 embeds  in  {\bf SLCG} as a reflective subcategory .

For a convex space $(X,\mathcal{C})$, Jin \cite{L.L16} defined a stratified $L$-cvs $(X,\omega_L(\mathcal{C}))$, where $\omega_L(\mathcal{C})$ is constructed from $\{\hat{a}\wedge \top_C|a\in L, C\in \mathcal{C}\}$ serving as a base.

\begin{proposition} \label{prpocpem} \cite{L.L16} $(X,\mathcal{C})\stackrel{f}\longrightarrow (Y,\mathcal{G})$ is a CP iff $(X,\omega_L(\mathcal{C}))\stackrel{f}\longrightarrow (Y,\omega_L(\mathcal{G}))$ is an LCP.
\end{proposition}

 \begin{proposition}  \label{thmwpi}  Given a  $U$-source $(X\stackrel{f_j}{\longrightarrow} U(X_j,\mathcal{C}_j))_{j\in J}$ in {\bf CVS} and its initial structure $\mathcal{C}$, then $\omega_L(\mathcal{C})$ is the initial structure of  the $U$-source $\big(X\stackrel{f_j}{\longrightarrow} (X_j,\omega_L(\mathcal{C}_j))\big)_{j\in J}$ in {\bf SLCVS}.
\end{proposition}

\begin{proof} We check that $\omega_L(\mathcal{C})$ fulfills the  two initial lift conditions in Definition \ref{defcvin}.

(1) Every $(X, \omega_L(\mathcal{C}))\stackrel{f_j}{\longrightarrow} (X_j, \omega_L(\mathcal{C}_j))$ is an LCP. It is concluded  from Proposition \ref{prpocpem}.

(2) For any mapping $(Y, \mathbb{G})\stackrel{f}{\longrightarrow} (X, \omega_L(\mathcal{C}))$ in {\bf SLCVS}, $f$ is an LCP iff each $f_j\circ f$ is an LCP.

The necessity is clear, only sufficiency remains for proving.  Let all $f_j\circ f$ be LCPs, we check below that $f$ is an LCP.

According to Proposition \ref{prpocp} and the fact that $\omega_L(\mathcal{C})$ has a base $\{\hat{a}\wedge \top_C|a\in L, C\in \mathcal{C}\}$, one require to verify that for all $a\in L$ and $C\in \mathcal{C}$, $f^\leftarrow(\hat{a}\wedge \top_C)\in \mathbb{G}$. Since $\mathcal{C}$ is generated by the family $\{f_j^\leftarrow(C_j)|j\in J, C_j\in \mathcal{C}_j\}$ as a subbase, we know that there exists $$\Big\{G_{ki}\in \mathcal{C}_{j_{ki}}|k\in K, i\in I_k, j_{ki}\in J \Big\}\ {\rm s.t.}\ C={\bigcup_{k\in K}}^\uparrow\bigcap_{i\in I_k}f_{j_{ki}}^\leftarrow(G_{ki}).$$
Furthermore, because all $f_{j_{ki}}\circ f$ are LCPs, it follows that $\top_{(f_{j_{ki}}\circ f)^\leftarrow (G_{ki})}\in \mathbb{G}$, and thus
$$f^\leftarrow(\hat{a}\wedge \top_C)=\hat{a}\wedge {\bigvee_{k\in K}}^\uparrow\bigwedge_{i\in I_k}\top_{(f_{j_{ki}}\circ f)^\leftarrow (G_{ki})}\in \mathbb{G}.$$Hence, $f$ is an LCP.
\end{proof}

Since the product structure is a special initial structure, we obtain the following corollary.
\begin{corollary} \label{lempp}  For a convex space $(X,\mathcal{C})$, $\omega_L(\mathcal{C}\times
\mathcal{C})=\omega_L(\mathcal{C})\times \omega_L(\mathcal{C})$.
\end{corollary}

\begin{proposition} \label{thmgg} Given  $(X,\cdot, \mathcal{C})\in $ {\bf CG}, then $(X,\cdot, \omega_L(\mathcal{C}))\in $ {\bf SLCG}.
\end{proposition}

\begin{proof} Since $(X,\cdot, \mathcal{C})\in $ {\bf CG}, then $k:(X\times X, \mathcal{C}\times\mathcal{C})\longrightarrow (X,\mathcal{C})$ is an LCP. According to Proposition \ref{prpocpem} and Corollary \ref{lempp} one get that
$$\Big(X\times X, \omega_L (\mathcal{C})\times
\omega_L (\mathcal{C})=\omega_L (\mathcal{C}\times
\mathcal{C})\Big)\stackrel{k}{ \longrightarrow} (X, \omega_L (\mathcal{C}))$$ to be LCP. Hence $(X,\cdot, \omega_L(\mathcal{C}))\in $ {\bf SLCG}.\end{proof}

Clearly,  the
correspondence $(X,\cdot, \mathcal{C})\mapsto (X,\cdot, \omega_L(\mathcal{C}))$ establishes an embedding functor, denoted as $\omega_L:{\bf CG}\longrightarrow
{\bf SLCG}$. From Propositions \ref{thmwpi} and \ref{thmgg} one deduce the next conclusion.

\begin{proposition} \label{thmwpig}  Given a $U$-source $\big((X,\cdot)\stackrel{f_j}{\longrightarrow} (X_j,\cdot,\mathcal{C}_j)\big)_{j\in J}$ in {\bf CG} and its initial structure  $\mathcal{C}$,  $\omega_L(\mathcal{C})$ is  the initial structure of the $U$-source $\big((X,\cdot)\stackrel{f_j}{\longrightarrow} (X_j,\cdot, \omega_L(\mathcal{C}_j))\big)_{j\in J}$ in {\bf SLCG}. We say that $\omega_L$ preserves the initial structures.
\end{proposition}

Since the join structures are special initial structures, we get the following corollary.

Let $CG(X,\cdot)$ represent  the collection of all convex groups over $(X,\cdot)$.

\begin{corollary}\label{corvv1}  For any $\{\mathcal{C}_k\}_{k\in K}\subseteq CG(X,\cdot)$, $\omega_L(\bigsqcup\limits_{k}\mathcal{C}_k)=\bigsqcup\limits_{k}\omega_L(\mathcal{C}_k)$.

\end{corollary}


Given $(X,\cdot, \mathbb{C})\in$ {\bf SLCG}, we define
$\rho_L(\mathbb{C})=\bigsqcup\{\mathcal{C}\in CG(X,\cdot)|\omega_L(\mathcal{C})\subseteq \mathbb{C}\}$. Then $(X,\cdot, \rho_L(\mathbb{C}))\in$ {\bf CG}.

\begin{proposition} If $(X,\cdot, \mathbb{C})\stackrel{f}\longrightarrow (Y,\cdot, \mathbb{G})$ is an LCP group homomorphism, then  $(X,\cdot,\rho_L(\mathbb{C}))$ $\stackrel{f}\longrightarrow (Y,\cdot,\rho_L(\mathbb{G}))$ is a CP group homomorphism.
\end{proposition}

\begin{proof} For any $\mathcal{G}\in CG(Y,\cdot)$ with $\omega_L(\mathcal{G})\subseteq \mathbb{G}$, let $f^\leftarrow(\mathcal{G})=\{f^\leftarrow(G)|G\in \mathcal{G}\}$. Obviously,  $(X,\cdot, f^\leftarrow(\mathcal{G}))$ forms a convex group. Since $(X,\cdot, \mathbb{C})\stackrel{f}\longrightarrow (Y,\cdot, \mathbb{G})$ is an LCP group homomorphism, $\omega_L(f^\leftarrow(\mathcal{G}))\subseteq \mathbb{C}$. Hence, $$f^\leftarrow(\mathcal{G})\subseteq \bigsqcup\{\mathcal{C}\in CG(X,\cdot)|\omega_L(\mathcal{C})\subseteq \mathbb{C}\}=\rho_L(\mathbb{C}).$$
It follows by  $\rho_L(\mathbb{G})=\bigsqcup\{\mathcal{G}\in CG(Y,\cdot)|\omega_L(\mathcal{G})\subseteq \mathbb{G}\}$ and Proposition \ref{prpocp}, we obtain that $(X,\cdot,\rho_L(\mathbb{C}))\stackrel{f}\longrightarrow (Y,\cdot,\rho_L(\mathbb{G}))$ is a CP group homomorphism.
\end{proof}

Therefore, the
correspondence $(X,\cdot,\mathbb{C})\mapsto (X,\cdot, \rho_L(\mathbb{C}))$ establishes a
concrete functor $\rho_L:{\bf SLCG}\longrightarrow
{\bf CG}$.

\begin{theorem} \label{thmref} $\rho_L$ is both a left inverse and a left adjoint of $\omega_L$. Hence, {\bf CG} can be reflectively embedded into {\bf SLCG} as a subcategory.
\end{theorem}

\begin{proof} Let $(X,\cdot,\mathcal{C})\in$ {\bf CG} and $(X,\cdot, \mathbb{C})\in$ {\bf
SLCG}.

(1) Obviously,  $\omega_L(\mathcal{G})\subseteq \omega_L(\mathcal{C})\Leftrightarrow \mathcal{G}\subseteq \mathcal{C}$ because $\omega_L$ is an embedding. Therefore $$\rho_L\circ \omega_L(\mathcal{C})=\bigsqcup\{\mathcal{G}\in CG(X,\cdot)|\omega_L(\mathcal{G})\subseteq  \omega_L(\mathcal{C})\}=\bigsqcup\{\mathcal{G}\in CG(X,\cdot)|\mathcal{G}\subseteq \mathcal{C}\}=\mathcal{C}.$$

(2) $\omega_L\circ \rho_L(\mathbb{C})\subseteq \mathbb{C}$. Indeed, \begin{eqnarray*} \omega_L\circ\rho_L(\mathbb{C})&=&\omega_L\Big(\bigsqcup\Big\{\mathcal{G}|\mathcal{G}\in CG(X,\cdot), \omega_L(\mathcal{G})\subseteq \mathbb{C}\Big\}\Big),\ {\rm by\ Corollary\ \ref{corvv1}}\\
&=&\bigsqcup\Big\{\omega_L(\mathcal{G})|\mathcal{G}\in CG(X,\cdot), \omega_L(\mathcal{G})\subseteq\mathbb{C}\Big\}\subseteq\mathbb{C}. \qedhere\end{eqnarray*}\end{proof}

\begin{remark} \label{remcate} In \cite{L.L16, PB19}, following  the functors between topological spaces and stratified $L$-topological spaces, Jin and Pang defined a pair of functors $\omega_L$ and $\rho_L$, between  convex spaces and stratified
$L$-convex spaces. However, the functor $\rho_L$ cannot map a stratified
$L$-convex group to a convex group, making it unsuitable for use in this paper. Thus, we propose a new approach, using the join structure, to redefine the functor $\rho_L$.
\end{remark}

\section{Concluding remarks}
This paper investigated the notion  of stratified $L$-convex groups and provided many equivalent characterizations. It has been proven that stratified $L$-convex groups and LCP group homomorphisms form a concrete category, denoted as {\bf SLCG}. This category has both initial and final structures, and consequently possesses substructures, join structures, product structures, and quotient structures. Furthermore, {\bf SLCG} has been shown to be a topological category over {\bf Grp}, and {\bf CG} can be embedded within {\bf SLCG} as a reflective subcategory.  To achieve this results, we employ different approaches and techniques compared to those used for (lattice-valued) topological groups and (lattice-valued) convex spaces.

The construction of uniformization and power objects is an important topic in topological groups, and we will conduct further research on {\bf SLCG} in this regard. In addition to the stratified $L$-convex structure, there are many other lattice-valued convex structures. We plan to investigate additional categories of lattice-valued convex groups and explore the relationships between them.

\end{document}